\documentclass{ifacconf}

\usepackage[english]{babel}
\usepackage{graphicx}      
\usepackage{natbib}        
\usepackage{subfigure}
\usepackage{color}
\usepackage{amssymb,amsmath,color}
\usepackage{srcltx,euscript}

\newcommand{\gs}{\geqslant}
\newcommand{\ilt}{\int\limits}

\begin{document}
\begin{frontmatter}
\title{Data-driven bandwidth choice for gamma kernel estimates of density derivatives
 on the positive semi-axis}
\author[Dobrovidov]{Alexander V. Dobrovidov,}
\author[Markovich]{Liubov A. Markovich}
\address[Dobrovidov]{Institute of Control Sciences, Russian Academy of Sciences,
Moscow, Russia (e-mail: dobrovid@ipu.ru).}
\address[Markovich]{Institute of Control Sciences, Russian Academy of Sciences,
Moscow, Russia (e-mail: kimo1@mail.ru).}
\begin{abstract}:
In some applications it is necessary to estimate derivatives of
probability densities defined on the positive semi-axis. The quality
of nonparametric estimates of the probability densities and their
derivatives are strongly influenced by smoothing parameters
(bandwidths). In this paper an expression for the optimal smoothing
parameter of the gamma kernel estimate of the density derivative is
obtained. For this parameter data-driven estimates based on methods
called "rule of thumb" are constructed. The
quality of the estimates is verified and demonstrated on examples of
density derivatives generated by Maxwell and Weibull distributions.
\end{abstract}
\end{frontmatter}

\section{Introduction}
In the paper of \cite{Chen:20} a nonparametric estimator of the
probability density function (pdf) for non-negative random variables
was obtained. The evaluation was built using  gamma kernels which
are asymmetric. In some applications like financial mathematics it
is necessary to estimate derivatives of the density. In this paper
the  estimate of the density derivative is constructed as a
derivative of the gamma kernel estimate.  In a previous paper of the
authors \cite{DobrovidovMarkovich:13} the statistical properties of
these kernel estimates defined on the positive semi-axis were
investigated. The asymptotic bias, variance and the integrated mean
squared error (MISE) as $n\to \infty$ were found. The smoothing
parameter that minimizes MISE, is called the optimal smoothing
parameter. It depends on functionals of the unknown true probability
density and its derivatives. Therefore, it is impossible to
calculate the exact value of this optimal smoothing parameter.
However,  the consistent estimates of this parameter can be
constructed using realizations of the underlying random variable
(r.v.) $X$. In the paper these estimates for the density derivative
are obtained by known methods called  "the rule of thumb" (see
\cite{Turlach}).

\section{Basic definitions and preliminary results}
Let a positive r.v. $X$ be described by the pdf $f(x)$, $x>0$. If
$f(x)$  is unknown, then the standard statistical problem is to
evaluate this density by an independent sample
$X^n=\{X_1,\ldots,X_n\}$. For the pdf defined on the positive
semi-axis in the paper of \cite{Chen:20} it is proposed a gamma
kernel estimator
 \vspace{-5mm}
\begin{eqnarray}\label{1}
\hat{f}(x) = \frac{1}{n}\sum_{i=1}^{n}K_{\rho_b(x),b}(X_i),
\end{eqnarray}
where \vspace{-8mm}
\begin{eqnarray*}K_{\rho_b(x),b}(t)=\frac{t^{\rho_b(x)-1}\exp(-t/b)}{b^{\rho_b(x)}\Gamma(\rho_b(x))}
\end{eqnarray*}
is a gamma kernel. Here $b$ is a smoothing parameter (bandwidth)
such that $b\rightarrow 0$ as $n\rightarrow\infty$, $\Gamma(\cdot)$
is the standard gamma function and
  \vspace{-4mm}
\begin{eqnarray*}
\rho_b(x)&=& \left\{
\begin{array}{ll}
\rho_1(x) = x/b, &   \mbox{if}\qquad x\geq 2b,
\\
\rho_2(x) =\left(x/2b\right)^2+1, & \mbox{if}\qquad x\in [0,2b).
\end{array}
\right.
\end{eqnarray*}
The support of the gamma kernel matches the support of the pdf to be
estimated. For convenience let us introduce two kernel functions on
the adjacent domains
   \vspace{-4mm}
\begin{eqnarray*}
 K_{\rho_1(x),b}(t)&=\frac{t^{x/b-1}\exp(-t/b)}{b^{x/b}\Gamma(x/b)},
  \quad    &\mbox{if}\qquad x\geq 2b,\\
 K_{\rho_2(x),b}(t)&=\frac{t^{(\frac{x}{2b})^2}\exp(-t/b)}{b^{(\frac{x}{2b})^2+1}\Gamma((\frac{x}{2b})^2+1)},
 \quad &\mbox{if}\qquad x\in [0,2b).
\end{eqnarray*}
The estimate of ${f}'(x)$ is constructed as the derivative of
$\hat{f}(x).$ Thus, we can write it as follows
 \vspace{1mm}
\begin{gather}\label{2}
\hat{f}'(x) = \frac{1}{n}\sum_{i=1}^{n}K'_{\rho_b(x),b}(X_i),
\end{gather}
where
 \vspace{-6mm}
\begin{eqnarray*}
&&K'_{\rho_b(x),b}(t)= \\
&=&\left\{
\begin{array}{ll}
K'_{\rho_1(x),b}(t)=\frac{1}{b}K_{\rho_1(x),b}(t)L_1,
& \mbox{if}\quad x\geq 2b,\\
K'_{\rho_2(x),b}(t)=\frac{x}{2b^2}K_{\rho_2(x),b}(t)L_2,&
\mbox{if}\quad x\in [0,2b),
\end{array}
\right.
\end{eqnarray*}
with
 \vspace{-6mm}
\begin{eqnarray*}
L_1=L_1(t;x)= \ln t - \ln b - \Psi(\rho_1(x)),\\
L_2=L_2(t;x)= \ln t - \ln b - \Psi(\rho_2(x)).
\end{eqnarray*}
Here $\Psi(x)$  denotes the Digamma function (logarithm of the
derivative of the Gamma-function). The asymptotic properties of the
density derivative estimate are determined by the following two
lemmas.
\par {\bf Lemma 1.}(expectation) {\it If \ $b\rightarrow 0$ then the
leading term of the
 mathematical expectation expansion for the density derivative
estimate  \eqref{2} equals \vspace{-5mm}
\begin{eqnarray}
\mathsf E(\hat{f'}(x)) = \left\{
\begin{array}{ll}
\mathsf E K'_{\rho_1(x),b}(X_1), & \mbox{if}\quad x\gs 2b, \notag\\
\mathsf E K'_{\rho_2(x),b}(X_1), & \mbox{if}\quad x\in [0,2b),\notag
\end{array}
\right.
\end{eqnarray}
where }
\begin{gather*}
\begin{array}{ll}
&\mathsf E K'_{\rho_1(x),b}(X_1)=(1/b)\mathsf E
K_{\rho_1(x),b}(X_1)L_1(X_1;x)\\
&= f'(x) + b\left(\frac{1}{12x^2} f(x) +
\frac{1}{4}f''(x)\right)+o(b),\\
&\mathsf E K'_{\rho_2(x),b}(X_1)=(x/2b^2)\mathsf E K_{\rho_2(x),b}(X_1)L_2(X_1;x)\\
&=f'(x)\left(\frac{x}{2b}-\frac{b}{6x}\right)+f''(x)\left(\frac{7x}{48}+
\frac{x^2}{2b}\right)+o(b).
\end{array}
\end{gather*}
Note that
under fixed $b$  the estimate   $\hat{f'}(x)$ in the small area
$x\in [0,2b)$ near zero has a  bias, which  grows as $x\rightarrow
0$ . However, in the asymptotic case when $b\rightarrow 0$ the right
boundary of this area $x=2b$ decreases also to zero. Therefore, it
is interesting to know the bias limit when $x$ and $b$ converge to
zero simultaneously, i.e. when ratio $x/b$ tends to some constant
$\kappa$ when $b\rightarrow 0$. Then the second expectation of the
estimate  will differ very small from the true density derivative.
The leading term of bias expansion may be written as
  \vspace{-5mm}
\begin{eqnarray*}
Bias(\hat{f}'(x))&=&b\left(\frac{f(x)}{12x^2}+
\frac{f''(x)}{4}\right)+ o(b), \quad i\!f \;x/b\rightarrow \infty,\\
Bias(\hat{f}'(x))&=&f'(x)\left(\frac{3\kappa^2-6\kappa-1}{6\kappa}\right)\\
&+& b f''(x)\left(\frac{7\kappa}{48}+
\frac{\kappa^2}{2}\right)+o(b), \quad i\!f \;x/b\rightarrow \kappa.
\end{eqnarray*}
If $x=2b$ then $\kappa=2$ and the estimate  bias in the right
boundary of the small area near zero will differ from true density
derivative in $(1/12) f'(2b)$.

As a global performance of the density derivative estimate \eqref{2}
we select the mean integrated squared error ($MISE$), which, as
known, equals
  \vspace{-5mm}
\begin{eqnarray}\label{3}
MISE(\hat{f}'(x))&=&
\mathsf E\int\limits_0^\infty(f'(x)-\hat{f}'(x))^2dx \\
&=& \int\limits_0^\infty \big[Bias^2(\hat{f}'(x)) +
Var(\hat{f}'(x))\big]dx. \notag
\end{eqnarray}
As the small area near zero will diminish with $n\rightarrow\infty,$
then the integral contribution in $MISE$ of the second part of the
bias in small area near zero will be negligible. Therefore, the only
integral squared bias of the main area of support is important.
Here it is
 \vspace{-4mm}
\begin{eqnarray*}
 IBias^2(\hat{f'}(x))
&=& \frac{b^2}{16} \ilt_{0}^\infty
\left(\frac{f(x)}{3x^2}+f''(x)\right)^2dx+o(b^2).
\end{eqnarray*}

Let us  proceed to calculate the variance of the derivative estimate.

{\bf Lemma 2.}(variance) {\it If  \ $b\rightarrow 0$ and $
nb^{3/2}\rightarrow \infty,$ then the leading term of
 variance expansion for density derivative estimate
\eqref{2} equals to}
 \vspace{-4mm}
 \begin{eqnarray*}&& Var(\hat{f}'(x))=
\\
&=& \frac{n^{-1}b^{-3/2}x^{-1/2}}{2\sqrt{\pi}}\left(\frac{f(x)}{2x}+b\left(\frac{f(x)}{4x^2}-\frac{f'(x)}{4x}\right)\right)+o(b)
\end{eqnarray*}

{\bf Theorem} ($MISE$). {\it If $b\rightarrow 0$ and \ $
nb^{3/2}\rightarrow \infty,$ integrals
  \vspace{-7mm}
\begin{eqnarray*}
\int_{0}^\infty \!\!\!
\left(\frac{f(x)}{3x^2}+f''(x)\right)^2\!\!\!dx, \int_{0}^\infty
\!\!\!\!\! x^{-3/2} f(x)dx
\end{eqnarray*}
are finite and $\int_{0}^\infty P(x) dx\neq0$, then the leading term of a MISE expansion for the
density derivative estimate  $\hat{f}'(x)$ equals to
 \vspace{-5mm}
\begin{eqnarray} \label{5}
 &&MISE(\hat f'(x))=\frac{b^2}{16}\int_{0}^\infty
 \left(\frac{f(x)}{3x^2}+f''(x)\right)^2dx \nonumber
\\
&+&\!\!  \int_0^\infty \frac{n^{-1}b^{-3/2}x^{-3/2}}{4\sqrt{\pi}}\left(f(x)+\frac{b}{2}\left(\frac{f(x)}{x}-f'(x)\right)\right)dx
\\
&+&o(b^2 + n^{-1}b^{-3/2}). \nonumber
\end{eqnarray}
Minimization of \eqref{5} in $b$ leads to a global optimal bandwidth
 \vspace{-5mm}
\begin{eqnarray}\label{6}
 b_0 = \left(\frac{3\int_0^\infty x^{-3/2}f(x)dx}{\sqrt{\pi}\int_{0}^\infty
(\frac{f(x)}{3x^2}+f''(x))^2dx}\right)^{2/7}n^{-2/7},
\end{eqnarray}
whose substitution into \eqref{5} results to an optimal {$MISE$} }
The restrictions on the integrals in the Theorem are fulfilled, for
example, for the family of $\chi^2$-distributions with a number of
degrees of freedom $m\gs 3.$ For $m=3$ we receive Maxwell
distribution, which will be investigated as true distribution in
simulation below. Proofs of the assertions are presented in
\cite{DobrovidovMarkovich:13}.

From expression for $MISE_{opt}$ it follows that nonparametric
estimate  \eqref{2} converges in mean square to true density
derivative with the rate $O(n^{-4/7}).$  This rate is certainly less
than the rate of convergence for the density $O(n^{-4/5})$, because the estimation of derivatives is more complex
than the estimation of the densities. A similar decrease in the rate of convergence
for the derivatives compared with the densities was observed in the
use of Gaussian kernel functions on the whole line.

\section{DATA-DRIVEN BANDWIDTH CHOICE}
The optimal smoothing parameter for the density derivative estimate,
defined by the formula \eqref{6}, depends on the unknown true
density $f(x)$ and its derivatives. Therefore, it is impossible to
calculate the true value of this parameter. However, one can
construct a non-parametric estimate of this parameter. Quite a lot
of methods for bandwidth estimation from the sample $X_1^n$ are
known. The  simplest and most convenient (in authors' opinion) are
methods called the rule of thumb (see \cite{Turlach}) and
cross-validation (see \cite{Rudeno:82}, \cite{Bowman:84}). The first
one is a parametric method. It gives a rough estimate by using in
\eqref{6} instead of the unknown true density a so-called reference
function, i.e. a density in the form of the kernel function.
Parameters of the latter can be set by any classical method of
parameter estimation. In this paper we use the method of moments.

The second approach for bandwidth estimation is to represent the
integrals in \eqref{6} in a form of expectation of some function
$\varphi(\cdot)$, i.e. as $\mathsf E\varphi(X)$. The function
$\varphi(\cdot)$ may depend also on the unknown true density $f(x)$
and its derivatives. We shall mark it  as
$\varphi\left(x,f^{(\alpha)}(x)\right),\; \alpha = 0,1,...$~. The
expectation of $\varphi(\cdot)$ equals to
 \vspace{-6mm}
\begin{eqnarray*}
\mathsf E\varphi\left(X,f^{(\alpha)}(X)\right)= \int\limits_0^\infty
\varphi\left(x,f^{(\alpha)}(x)\right)f(x)dx
\end{eqnarray*}
 by definition and it can be approximated by the arithmetic mean
 \vspace{-7mm}
\begin{eqnarray}\label{7}
\mathsf E\varphi\left(X,f^{(\alpha)}(X)\right)\approx
\frac{1}{n}\sum\limits_{i=1}^{n}\varphi\left(X_i,f^{(\alpha)}(X_i)\right),
\end{eqnarray}
converging to it as $n\to \infty$. The unknown density and its
derivatives in \eqref{7} are replaced by their gamma kernel
estimates, but in the form of Cross-validation
 \vspace{-6mm}
\begin{eqnarray}\label{8}
\hat{f}_i^{(\alpha)}(X_i) =
\frac{1}{n-1}\sum\limits_{\genfrac{}{}{0pt}{}{j=1}{j\neq i}}^n
K_{\rho(X_i)}^{(\alpha)}(X_j), \quad \alpha = 0,1,...~.
\end{eqnarray}
These estimates will also depend on the unknown smoothing parameters
which are substituted by rough estimates from  the rule of thumb.
Experience has shown that roughness of bandwidth estimation on the
second level does not  affect too much on the accuracy of estimation
of densities and their derivatives, \cite{Devroye:85}.

Let us  find the value of the smoothing parameter following the rule
of thumb. For this suppose we choose the gamma density
 \vspace{-8mm}
\begin{eqnarray} \label{9a}
f(x)=\frac{x^{\rho-1}\exp(-\frac{x}{b})}{b^{\rho}\Gamma(\rho)}
\end{eqnarray}
as a reference function.  The first moment and the variance of it
are $\rho b$ and $\rho b^2$, respectively. According to the method
of moments, we have to equate them to the first sample moment $\bar
m = n^{-1}\sum_{i=1}^n X_i $ and the sample variance $\bar D =n^{-1}
\sum_{i=1}^{n}(X_i - \bar m)^2$, correspondingly. Then we obtain for
the parameters of \eqref{9a} following simple expressions
 \vspace{-7mm}
\begin{eqnarray}\label{10}
 b_m = \bar D/\bar m,\quad \rho_m = (\bar m)^2/\bar D.
\end{eqnarray}
 Fig. 1  gives an idea of the relative closeness of the true (Maxwell) and references
 (gamma) densities, when the first two moments of the corresponding distributions
  are almost identical.

 \begin{figure}
\begin{center}
\includegraphics[width=8.4cm]{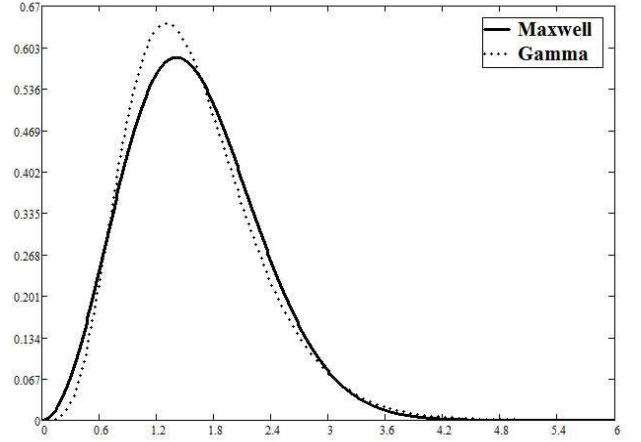}
\caption{Gamma density with parameters $b_m,\rho_m$ and the true
Maxwell density functions} \label{fig:1}
\end{center}
\end{figure}

Now one can replace the unknown density in \eqref{6} by a gamma
density \eqref{9a} with known parameters \eqref{10}. Then the
numerator of $b_0$ is
 \vspace{-6mm}
\begin{eqnarray}&&
I_{n}=\! \int\limits_0^\infty\! \!\frac{3}{\sqrt{\pi}}x^{-3/2}f(x)dx \\ \nonumber
&=&\frac{3}{\sqrt{\pi} b_{m}^{\rho_{m}}\Gamma(\rho_{m})}\ilt_0^\infty t^{\rho_{m}-\frac{5}{2}}\exp\left(\frac{-t}{b_{m}}\right)dt
\end{eqnarray}
The denominator of \eqref{6} equals
 \vspace{-6mm}
\begin{eqnarray}&&
I_{d}=\ilt_0^\infty \left(\frac{f(x)}{3x^2}+f''(x)\right)^2dx\\
\nonumber & =& \ilt_0^\infty
\left(\frac{t^{\rho_m-1}\exp(\frac{-t}{b_m})}{3t^2b_m^\rho\Gamma(\rho_m)}+\frac{d^2}{dt^2}\left(\frac{t^{\rho_m-1}\exp(\frac{-t}{b_m})}{b_m^\rho\Gamma(\rho_m)}\right)\right)^2dt\\
\nonumber
&=&\!\frac{\Gamma(\rho_m\!-\!\frac{5}{2})(b_m^{\rho_m}(4b^2\!-\!12\rho_m\!+\!48)-81\rho_m+27\rho_m^2+54)}{72\sqrt{\pi}\Gamma(\rho_m)(\rho_m-1)(\rho_m-2)b_m^5}.\nonumber
\end{eqnarray}
Thus, collecting  both terms in one formula, we obtain an expression
for the bandwidth from the rule of thumb as
 \vspace{-8mm}
\begin{eqnarray}\label{12}
b_{0G}=\left(\frac{I_n}{I_d}\right)^{2/7}\cdot n^{-2/7}.
\end{eqnarray}
\section{Simulation results}
In the simulation experiment we have selected a
 Maxwell density ($\sigma=1$)
  \vspace{-8mm}
\begin{eqnarray*}
f_M(x)=\frac{\sqrt{2}x^2\exp(-x^2/2\sigma^2)}{\sigma^3\sqrt{\pi}},
\end{eqnarray*}
with its derivative
  \vspace{-8mm}
\begin{eqnarray*}
f'_M(x)=-\frac{\sqrt{2}x\exp(-x^2/2\sigma^2)(x^2-2\sigma^2)}{\sigma^5\sqrt{\pi}}
\end{eqnarray*}
 and a Weibull density ($s=3$)
 \vspace{-7mm}
\begin{eqnarray*} f_W(x)=sx^{s-1}\exp(-x^s),
\end{eqnarray*}
with its derivative
  \vspace{-8mm}
\begin{eqnarray*}
f'_W(x)=-sx^{s-2}\exp(-x^s)(sx^s-s+1)
 \end{eqnarray*}
as the unknown true functions to be estimated. Parameters of the
densities are selected to satisfy the integrated conditions of the
Theorem above. We generated Maxwell and Weibull samples $X_1^n$ with
sample length $n=1000$ using standard generators (see Fig. 2).
\begin{figure}
\begin{center}
\includegraphics[width=8.4cm]{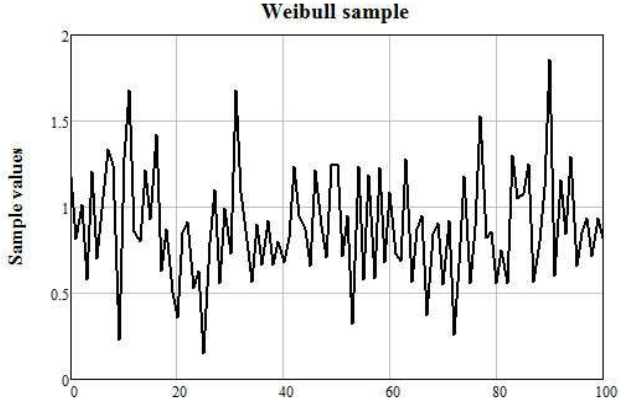}
\caption{Sample from a Weibull distribution} \label{fig:2}
\end{center}
\end{figure}

Based on these samples and the method of moments \eqref{10}, we
obtain the parameter values  $b_{m}$ and $\rho_{m}$ for
the reference gamma density \eqref{9a}. Substitution of it into the
expression \eqref{6} for the optimal bandwidth $b_0$ leads to the
estimate by the rule of thumb
\begin{gather}\label{16}
\hat b_{0G}=I_n/I_d.
 \end{gather}
Similar computations were performed for the sample from a Maxwell
distribution. C

Density derivative estimates corresponding to calculated bandwidths
are represented on Fig.~ 3 and Fig. 4. Here
solid line corresponds to the true density derivative to be
estimated, dashed line corresponds to rule of thumb estimates.

\begin{figure}
\begin{center}
\includegraphics[width=8.4cm]{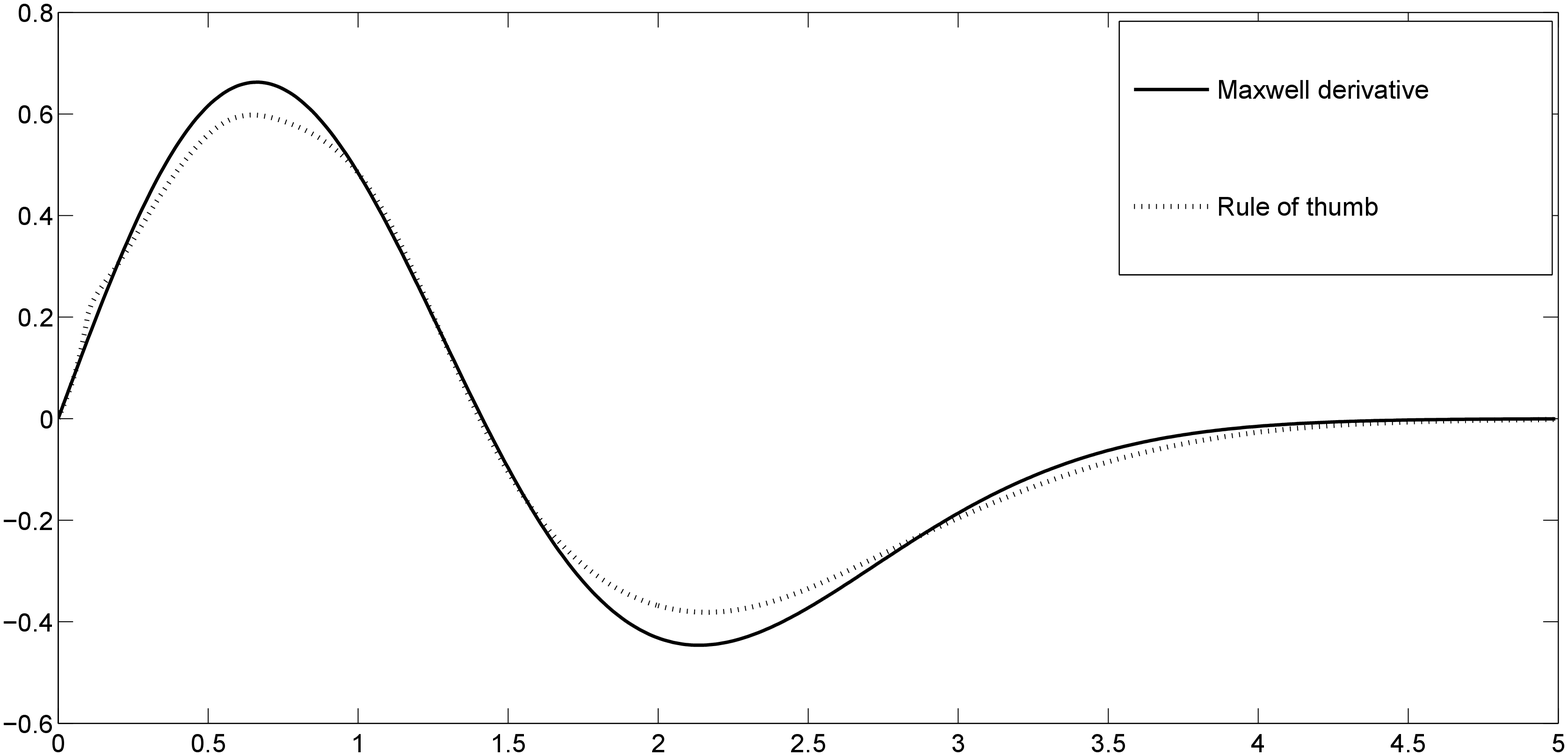}
\caption{Nonparametric estimates of the Maxwell density derivative
function for $n=10000$. The true $f'_M(x)$ (solid line), estimate
from the rule of thumb  (dashed line)} \label{fig:3}
\end{center}
\end{figure}
\begin{figure}
\begin{center}
\includegraphics[width=8.4cm]{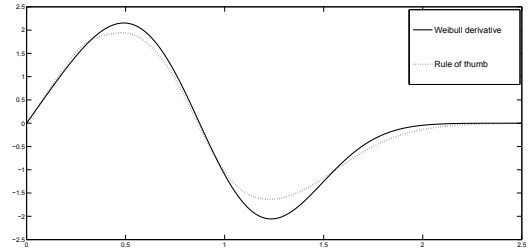}
\caption{Nonparametric estimates of the Weibull density derivative
function for $n=10000$. The true $f'_W(x)$ (solid line), estimate
from the rule of thumb (dashed  line)} \label{fig:4}
\end{center}
\end{figure}
Quantitatively, the estimation error is determined by the value
$\varkappa$, defined by the following formula:
\begin{gather}\label{15}
\varkappa=x_{step}\sum((f(x)-\widehat{f}(x))^2).
\end{gather}
The mean  error
\begin{gather*}
\overline{\varkappa}=M^{-1}\sum_{i=1}^{M}\varkappa_i
\end{gather*}
calculated over $M=100$ repeated experiments and standard deviations in
brackets are represented in Table 3 and Table 4.

\begin{table}[htbp]
    \begin{center}
    \caption{Error $\overline{\varkappa}$ for the Maxwell distribution} {} \label{tab:3}
    \vspace{-4mm}
    \end{center}
    \begin{center}
    \begin{tabular}{|c|c|c|}
    \hline m &rule of thumb  \\
    \hline $m=100$  &0.0082\\
    \hline$m=1000$  &0.02302\\
    \hline$m=2000$  &0.015314\\
    \hline
    \end{tabular}
    \end{center}
    \end{table}
\vspace{4mm}
\begin{table}[htbp]
    \begin{center}
    \caption{Error $\overline{\varkappa}$ for the Weibull distribution} {} \label{tab:4}
    \end{center}
\vspace{-4mm}
    \begin{center}
    \begin{tabular}{|c|c|}
     \hline m &rule of thumb  \\
    \hline $m=100$ &0.61945\\
    \hline$m=1000$ &0.25808 \\
    \hline$m=2000$  &0.17994 \\
    \hline
    \end{tabular}
    \end{center}
    \end{table}
From the tables and graphics one can see that from two methods of
bandwidth selection (the rule of thumb and cross-validation) better
results gives the latter one.

\section{CONCLUSIONS}
Estimation of the probability characteristics of the positive random
variables is required in the theory of signal processing, the
financial and actuarial mathematics and other important
applications. The positivity of the distribution support of the
observed random variables results in a significant complication of
the models compared to the case of an unbounded support.  We present
nonparametric kernel estimates of the  densities and their
derivatives based on the asymmetric gamma
 kernels. The  main part of the paper is devoted to the construction and evaluation of
the optimal  smoothing parameters (bandwidths) in the kernel
 density derivative estimates  by samples of the independent random variables.
 For this purpose two well-known method
  such as the rule of thumb is used. It should be noted
  that in the case of asymmetric
support even for the simple rule of thumb method the expression for
the smoothing parameter \eqref{12} becomes quite cumbersome in
comparison with the known corresponding expression based on the
Gaussian reference function defined on the whole line. The estimates of the kernel density derivatives will be
    used for the nonparametric estimation of logarithmic derivatives of the density
    determined on the positive real axis.   The latter can be used for the problems
    of unsupervised nonparametric signal processing.

\bibliography{reference}
\end{document}